\documentclass[10pt,twocolumn,twoside]{IEEEtran}%
\usepackage{amsmath}
\usepackage{graphicx}
\usepackage{amsfonts}
\usepackage{amssymb}%
\usepackage{epstopdf}%
\makeatletter

\@addtoreset{theorem}{section}

\@addtoreset{corollary}{section}

\@addtoreset{lemma}{section}

\@addtoreset{definition}{section}
\makeatother

\begin{document}

\title{Energy-aware Vehicle Routing in Networks with Charging
Nodes\thanks{The authors' work is supported in part by NSF under
Grant CNS-1139021, by AFOSR under grant FA9550-12-1-0113, by
ONR under grant N00014-09-1-1051, and by ARO under Grant W911NF-11-1-0227.}}
\author{\textbf{ C.~G. Cassandras, T. Wang }and\textbf{ S. Pourazarm}\\Division of Systems Engineering and\\Center for Information and Systems Engineering, Boston University\\\texttt{cgc@bu.edu, renowang0823@gmail.com, sepid@bu.edu}}
\maketitle

\begin{abstract}                
We study the problem of routing vehicles with energy constraints through a network where there are at least some charging nodes.
We seek to minimize the total elapsed time for vehicles to reach their destinations by determining
routes as well as recharging amounts when the vehicles do not have adequate
energy for the entire journey. For a single vehicle, we formulate a
mixed-integer nonlinear programming (MINLP) problem and derive properties of
the optimal solution allowing it to be decomposed into two simpler problems.
For a multi-vehicle problem, where traffic congestion effects are included, we
use a similar approach by grouping vehicles into \textquotedblleft
subflows.\textquotedblright\ We also provide an alternative flow optimization
formulation leading to a computationally
simpler problem solution with minimal loss in accuracy. Numerical results are included to illustrate these approaches.
\end{abstract}

\section{Introduction}

\label{sec1}

The increasing presence of Battery-Powered Vehicles (BPVs), such as Electric
Vehicles (EVs), mobile robots and sensors, has given rise to novel issues in
classical network routing problems [\cite{Laporte91}]. More generally, when the entities in the
network are characterized by physical attributes exhibiting a dynamic
behavior, this behavior can play an important role in the routing decisions.
In the case of BPVs, the physical attribute is \emph{energy} and there are
four BPV characteristics which are crucial in routing problems: limited
cruising range, long charge times, sparse coverage of charging stations, and
the BPV energy recuperation ability [\cite{Artmeier2010}] which can be
exploited. In recent years, the vehicle routing literature has been enriched
by work aiming to accommodate these BPV characteristics. For example, by
incorporating the recuperation ability of EVs (which leads to negative energy
consumption on some paths), extensions to general shortest-path algorithms are
proposed in \cite{Artmeier2010} that address the energy-optimal routing
problem. The energy requirements in this problem are modeled as constraints
and the proposed algorithms are evaluated in a prototypical navigation system.
Extensions provided in \cite{Eisner2011} employ a generalization of Johnson's
potential shifting technique to make Dijkstra's algorithm applicable to the
negative edge cost shortest-path problem so as to improve the results and
allow for route planning of EVs in large networks. This work, however, does
not consider the presence of charging stations, modeled as nodes in the
network. Charging times are incorporated into a multi-constrained optimal path
planning problem in \cite{Siddiqi2011}, which aims to minimize the length of
an EV's route and meet constraints on total traveling time, total time delay
due to signals, total recharging time and total recharging cost. A particle
swarm optimization algorithm is used to find a suboptimal solution. In this
formulation, however, recharging times are simply treated as parameters and
not as controllable variables. In \cite{Kuller2011}, algorithms for several
routing problems are proposed, including a single vehicle routing problem with
inhomogeneously priced refueling stations for which a dynamic programming
based algorithm is proposed to find a least cost path from source to
destination. More recently, an EV Routing Problem with Time Windows and
recharging stations (E-VRPTW) was proposed in \cite{TechrptMichael}, where an
EV's energy constraint is first introduced into vehicle routing problems and
recharging times depend on the battery charge of the vehicle upon arrival at
the station. Controlling recharging times is circumvented by simply forcing
vehicles to be always fully recharged. In the Unmanned Autonomous Vehicle
(UAV) literature, \cite{Sunder2012} consider a UAV routing problem with
refueling constraints. In this problem, given a set of targets and depots the
goal is to find an optimal path such that each target is visited by the UAV at
least once while the fuel constraint is never violated. A Mixed-Integer
Nonlinear Programming (MINLP) formulation is proposed with a heuristic
algorithm to determine feasible solutions.

In this paper, our objective is to investigate a vehicle total traveling time
minimization problem (including both the time on paths and at charging
stations), where an energy constraint is considered so that the vehicle is not
allowed to run out of power before reaching its destination. We view this as a
network routing problem where vehicles control not only their routes but also
times to recharge at various nodes in the network. Our contributions are
twofold. First, for the single energy-aware vehicle routing problem,
formulated as a MINLP, we show that there are properties of the optimal
solution and the energy dynamics allowing us to decompose the original problem
into two simpler problems with inhomogeneous prices at charging nodes but
homogeneous charging speeds. Thus, we separately determine route selection
through a Linear Programming (LP) problem and then recharging amounts through
another LP or simple optimal control problem. Since we do not impose full
recharging constraints, the solutions obtained are more general than, for
example, in \cite{TechrptMichael} and recover full recharging when this is
optimal. Second, we study a multi-vehicle energy-aware routing problem, where
a traffic flow model is used to incorporate congestion effects. This
system-wide optimization problem appears to have not yet attracted much
attention. By grouping vehicles into \textquotedblleft
subflows\textquotedblright\ we are once again able to decompose the problem
into route selection and recharging amount determination, although we can no
longer reduce the former problem to an LP. Moreover, we provide an alternative
flow-based formulation such that each subflow is not required to follow a
single end-to-end path, but may be split into an optimally determined set of
paths. This formulation reduces the computational complexity of the MINLP
problem by orders of magnitude with numerical results showing little or no
loss in optimality.

The structure of the paper is as follows. In Section \ref{sec2}, we introduce
and address the single-vehicle routing problem and identify properties which
lead to its decomposition. In Section \ref{sec3}, the multi-vehicle routing
problem is formulated, first as a MINLP and then as an alternative flow
optimization problem. Simulation examples are included for the multi-vehicle
routing problem illustrating our approach and providing insights on the
relationship between recharging speed and optimal routes. Finally, conclusions
and further research directions are outlined in Section \ref{sec4}.

\section{Single Vehicle Routing}

\label{sec2}

\begin{figure}[ptbh]
\begin{center}
\includegraphics[scale=0.5]{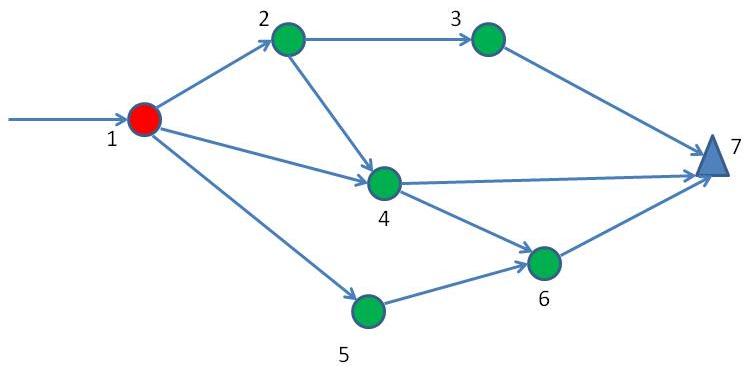}
\end{center}
\caption{A 7-node network example for routing with recharging nodes.}%
\label{SampleNet}%
\end{figure}

We assume that a network is defined as a directed graph $G=(\mathcal{N}%
,\mathcal{A})$ with $\mathcal{N}=\{1,\dots,n\}$ and $|\mathcal{A}|=m$ (see
Fig. \ref{SampleNet}). Node $i\in\mathcal{N}/\{n\}$ represents a charging
station and $(i,j)\in\mathcal{A}$ is an arc connecting node $i$ to $j$ (we
assume for simplicity that all nodes have a charging capability, although this
is not necessary). We also define $I(i)$ and $O(i)$ to be the set of start
nodes (respectively, end nodes) of arcs that are incoming to (respectively,
outgoing from) node $i$, that is, $I(i)=\{j\in\mathcal{N}|(j,i)\in
\mathcal{A}\}$ and $O(i)=\{j\in\mathcal{N}|(i,j)\in\mathcal{A}\}$.

We are first interested in a single-origin-single-destination vehicle routing
problem. Nodes 1 and $n$ respectively are defined to be the origin and
destination. For each arc $(i,j)\in\mathcal{A}$, there are two cost
parameters: the required traveling time $\tau_{ij}$ and the required energy
consumption $e_{ij}$ on this arc. Note that $\tau_{ij}>0$ (if nodes $i$ and
$j$ are not connected, then $\tau_{ij}=\infty$), whereas $e_{ij}$ is allowed
to be negative due to a BPV's potential energy recuperation effect
[\cite{Artmeier2010}]. Letting the vehicle's charge capacity be $B$, we assume
that $e_{ij}<B$ for all $(i,j)\in\mathcal{A}$. Since we are considering a
single vehicle's behavior, we assume that it will not affect the overall
network's traffic state, therefore, $\tau_{ij}$ and $e_{ij}$ are assumed to be
fixed depending on given traffic conditions at the time the single-vehicle
routing problem is solved. Clearly, this cannot apply to the multi-vehicle
case in the next section, where the decisions of multiple vehicle routes
affect traffic conditions, thus influencing traveling times and energy
consumption. Since the BPV has limited battery energy it may not be able to
reach the destination without recharging. Thus, recharging amounts at charging
nodes $i\in\mathcal{N}$ are also decision variables.

We denote the selection of arc $(i,j)$ and energy recharging amount at node
$i$ by $x_{ij}\in\{0,1\}$, $i,j\in\mathcal{N}$ and $r_{i}\geq0$,
$i\in\mathcal{N}/\{n\}$, respectively. Moreover, since we take into account
the vehicle's energy constraints, we use $E_{i}$ to represent the vehicle's
residual battery energy at node $i$. Then, for all $E_{j},\,j\in O(i)$, we
have:
\[
E_{j}=\left\{
\begin{array}
[c]{ll}%
E_{i}+r_{i}-e_{ij} & \text{if }x_{ij}=1\\
0 & \text{otherwise}%
\end{array}
\right.
\]
which can also be expressed as
\[
E_{j}=\sum_{i\in I(j)}(E_{i}+r_{i}-e_{ij})x_{ij},\quad x_{ij}\in\{0,1\}
\]
The problem objective is to determine a path from $1$ to $n$, as well as
recharging amounts, so as to minimize the total elapsed time for the vehicle
to reach the destination. Fig. \ref{SampleNet} is a sample network for this
vehicle routing problem. We formulate a MINLP problem as follows:
\begin{gather}
\min_{x_{ij},r_{i},\,\,i,j\in\mathcal{N}}\quad\sum_{i=1}^{n}\sum_{j=1}^{n}%
\tau_{ij}x_{ij}+\sum_{i=1}^{n}\sum_{j=1}^{n}r_{i}gx_{ij}\label{obj}\\
s.t.\quad\sum_{j\in O(i)}x_{ij}-\sum_{j\in I(i)}x_{ji}=b_{i},\quad\text{for
each }i\in\mathcal{N}\label{flowConv}\\
b_{1}=1,\,b_{n}=-1,\,b_{i}=0,\text{ for }i\neq1,n\label{bi}\\
E_{j}=\sum_{i\in I(j)}(E_{i}+r_{i}-e_{ij})x_{ij},\text{ for }j=2,\dots
,n\label{EiConv}\\
0\leq E_{i}\leq B,\quad E_{1}\text{ given},\text{ for each }i\in
\mathcal{N}\label{Ei}\\
x_{ij}\in\{0,1\},\quad r_{i}\geq0 \label{controls}%
\end{gather}
where $g$ is the charging time per energy unit, i.e., the reciprocal of a
fixed charging rate. The constraints (\ref{flowConv})-(\ref{bi}) stand for the
flow conservation [\cite{Bertsimas}], which implies that only one path
starting from node $i$ can be selected, i.e., $\sum_{j\in O(i)}x_{ij}\leq1$.
It is easy to check that this also implies $x_{ij}\leq1$ for all $i,j$ since
$b_{1}=1$, $I(1)=\varnothing$. Constraint (\ref{EiConv}) represents the
vehicle's energy dynamics where the only non-linearity in this formulation
appears. Finally, (\ref{Ei}) indicates that the vehicle cannot run out of
energy before reaching a node or exceed a given capacity $B$. All other
parameters are predetermined according to the network topology.

\subsection{Properties}

Rather than directly tackling the MINLP problem, we derive some key properties
which will enable us to simplify the solution procedure. The main difficulty
in this problem lies in the coupling of the decision variables, $x_{ij}$ and
$r_{i}$, in (\ref{EiConv}). The following lemma will enable us to exclude
$r_{i}$ from the objective function by showing that the difference between the
total recharging energy and the total energy consumption while traveling is
given only by the difference between the vehicle's residual energy at the
destination and at the origin.

\textbf{Lemma 1: }Given (\ref{obj})-(\ref{controls}),\textbf{ }
\begin{equation}
\sum_{i=1}^{n}\sum_{j=1}^{n}(r_{i}x_{ij}-e_{ij}x_{ij})=E_{n}-E_{1}
\label{Lemma1}%
\end{equation}
\emph{Proof}: From (\ref{EiConv}), we sum up both sides to get:
\begin{equation}
\sum_{j=2}^{n}E_{j}-\sum_{j=2}^{n}\sum_{i\in I(j)}E_{i}x_{ij}=\sum_{j=2}%
^{n}\sum_{i\in I(j)}(r_{i}-e_{ij})x_{ij} \label{EiCon1}%
\end{equation}
Moreover, we can write
\[
\sum_{j=2}^{n}\sum_{i\in I(j)}E_{i}x_{ij}=\sum_{i\in I(2)}E_{i}x_{i2}%
+\cdots+\sum_{i\in I(n)}E_{i}x_{in}%
\]
representing the sum of $E_{i}$ on the selected path from node $1$ to $n$,
excluding $E_{n}$. On the other hand, from (\ref{EiConv}) we have $E_{i}=0$
for any node $i$ not selected on the path. Therefore, $\sum_{j=2}^{n}E_{j}$ is
the sum of $E_{i}$ on the selected path from node $1$ to $n$, excluding
$E_{1}$. It follows that
\begin{equation}
\sum_{j=2}^{n}E_{j}-\sum_{j=2}^{n}\sum_{i\in I(j)}E_{i}x_{ij}=E_{n}-E_{1}
\label{En-E1}%
\end{equation}
Returning to (\ref{EiCon1}), we use (\ref{En-E1}) and observe that all terms
in the double sum $\sum_{i=1}^{n}\sum_{j=1}^{n}(r_{i}-e_{ij})x_{ij}$ are zero
except for those with $i\in I(j)$, we get
\begin{align}
&  \sum_{i=1}^{n}\sum_{j=1}^{n}(r_{i}-e_{ij})x_{ij}=\sum_{j=2}^{n}\sum_{i\in
I(j)}(r_{i}-e_{ij})x_{ij}\nonumber\\
&  =\sum_{j=2}^{n}E_{j}-\sum_{j=2}^{n}\sum_{i\in I(j)}E_{i}x_{ij}=E_{n}%
-E_{1}\nonumber
\end{align}
which proves the lemma.$\blacksquare$

In view of Lemma 1, we can replace $\sum_{i=1}^{n}\sum_{j=1}^{n}r_{i}gx_{ij}$
in (\ref{obj}) by $(E_{n}-E_{1})g+\sum_{i=1}^{n}\sum_{j=1}^{n}e_{ij}gx_{ij}$
and eliminate the presence of $r_{i}$, $i=2,\ldots,n-1$, from the objective
function. Note that $E_{1}$ is given, leaving us only with the task of
determining the value of $E_{n}$. Now, let us investigate the recharging
energy amounts $r_{i}^{\ast}$, $i=1,\ldots,n-1$, in an optimal policy. There
are two possible cases: $(i)$ $\sum_{i}r_{i}^{\ast}>0$, i.e., the vehicle has
to get recharged at least once, and $(ii)$ $\sum_{i}r_{i}^{\ast}=0$, i.e.,
$r_{i}^{\ast}=0$ for all $i$ and the vehicle has adequate energy to reach the
destination without recharging. For Case $(i)$, we establish the following lemma.

\textbf{Lemma 2: } If $\sum_{i}r_{i}^{\ast}>0$ in the optimal routing policy,
then $E_{n}^{\ast}=0$.\newline\emph{Proof}: We use a contradiction argument.
Assume we have already achieved an optimal route where $E_{n}^{\ast}>0$ and
the objective function is $J^{\ast}=\sum_{i\in P}(\tau_{i,i+1}+r_{i}^{\ast}g)$
for an optimal path denoted by $P$. Without loss of generality, we re-index
nodes so that we may write $P=\{1,\ldots,n\}$. Then, each $i\in P$ such that
$i<n$ on this optimal path satisfies:
\begin{equation}
E_{i+1}^{\ast}=E_{i}^{\ast}+r_{i}^{\ast}-e_{i,i+1} \label{EioneRoute}%
\end{equation}
\newline Consider first the case where $r_{n-1}^{\ast}>0$. Let us perturb the
current policy as follows: $r_{n-1}^{^{\prime}}=r_{n-1}^{\ast}-\Delta$, and
$r_{i}^{^{\prime}}=r_{i}^{\ast}$ for all $i<n-1$, where $\Delta>0$. Then, from
(\ref{EioneRoute}), we have
\[
E_{n}^{\ast}=E_{1}+\sum_{i=1}^{n-1}(r_{i}^{\ast}-e_{i,i+1})
\]
Under the perturbed policy,
\begin{align}
E_{n}^{^{\prime}}  &  =E_{1}+\sum_{i=1}^{n-1}(r_{i}^{^{\prime}}-e_{i,i+1}%
)\nonumber\\
&  =E_{1}+\sum_{i=1}^{n-1}(r_{i}^{\ast}-e_{i,i+1})-\Delta=E_{n}^{\ast}%
-\Delta\nonumber\\
E_{i}^{^{\prime}}  &  =E_{i}^{\ast},\text{ for all }i<n\nonumber
\end{align}
and, correspondingly,
\[
J^{^{\prime}}=\sum_{i=1}^{n-1}(\tau_{i,i+1}+r_{i}^{^{\prime}}\cdot
g)=\sum_{i=1}^{n-1}(\tau_{i,i+1}+r_{i}^{\ast}\cdot g)-\Delta g=J^{\ast}-\Delta
g
\]
Since $E_{n}^{\ast}>0$, we may select $\Delta>0$ sufficiently small so that
$E_{n}^{^{\prime}}>0$ and the perturbed policy is still feasible. However,
$J^{^{\prime}}=J^{\ast}-\Delta\cdot g<J^{\ast}$, which leads to a
contradiction to the assumption that the original path was optimal.

Next, consider the case where $r_{n-1}^{\ast}=0$. Then, due to $E_{n}^{\ast
}>0$ and $e_{i,i+1}>0$ for all $i\in P$, we can always find some $j\in P,$
$j<n$ such that $E_{j}^{\ast}>0$, $r_{j-1}^{\ast}>0$ and $r_{k}^{\ast}=0$ for
$k \geqslant j$. Thus, still due to (\ref{EioneRoute}), we have
\[
E_{j}^{\ast}=E_{n}^{\ast}+\sum_{k=j}^{n-1}e_{k,k+1}>0
\]
At this time, since $r_{j-1}^{\ast}>0$, the argument is similar to the case
$r_{n-1}^{\ast}>0$, leading again to the same contradiction argument and the
lemma is proved. $\blacksquare$

Turning our attention to Case $(ii)$ where $r_{i}^{\ast}=0$ for all
$i\in\{1,\ldots,n\}$, observe that the problem (\ref{obj}) can be transformed
to
\begin{gather}
\min_{_{x_{ij,\text{ }}i,j\in\mathcal{N}}}\quad\sum_{i=1}^{n}\sum_{j=1}%
^{n}\tau_{ij}x_{ij}\label{obj2}\\
s.t.\quad\sum_{j\in O(i)}x_{ij}-\sum_{j\in I(i)}x_{ji}=b_{i},\quad\text{for
each }i\in\mathcal{N}\nonumber\\
b_{1}=1,\,b_{n}=-1,\,b_{i}=0,\text{ for }i\neq1,n\nonumber\\
E_{j}=\sum_{i\in I(j)}(E_{i}-e_{ij})x_{ij},\text{ for }j=2,\dots
,n\label{EiConv2}\\
0\leq E_{i}\leq B,\quad E_{0}\text{ given},\text{ for each }i\in
\mathcal{N}\nonumber\\
x_{ij}\in\{0,1\} \label{controls2}%
\end{gather}
In this case, the constraint (\ref{EiConv2}) gives
\[
\sum_{j=2}^{n}E_{j}-\sum_{j=2}^{n}\sum_{i\in I(j)}E_{i}=-\sum_{j=2}^{n}%
\sum_{i\in I(j)}e_{ij}x_{ij}%
\]
Using (\ref{En-E1}) and $E_{i}\geq0$, we have
\[
E_{n}=E_{1}-\sum_{j=2}^{n}\sum_{i\in I(j)}e_{ij}x_{ij}\geq0
\]
and it follows that%
\begin{equation}
\sum_{i=1}^{n}\sum_{j=1}^{n}e_{ij}x_{ij}\leq E_{1} \label{leE1}%
\end{equation}
With (\ref{leE1}) in place of (\ref{EiConv2}), the determination of
$x_{ij}^{\ast}$ boils down to an integer \emph{linear} programming problem in
which only variables $x_{ij}$, $i,j\in\mathcal{N}$, are involved, a much
simpler problem.

We are normally interested in Case $(i)$, where some recharging decisions must
be made, so let us assume the vehicle's initial energy is not large enough to
reach the destination. Then, in view of Lemmas 1 and 2, we have the following theorem.

\textbf{Theorem 1: } If $\sum_{i}r_{i}^{\ast}>0$ in the optimal policy, then
$x_{ij}^{\ast}$, $i,j\in\mathcal{N}$, in the original problem (\ref{obj}) can
be determined by solving a linear programming problem:
\begin{gather}
\min_{_{x_{ij,\text{ }}i,j\in\mathcal{N}}}\quad\sum_{i=1}^{n}\sum_{j=1}%
^{n}(\tau_{ij}+e_{ij}g)x_{ij}\label{obj3}\\
s.t.\quad\sum_{j\in O(i)}x_{ij}-\sum_{j\in I(i)}x_{ji}=b_{i},\quad\text{for
each }i\in\mathcal{N}\nonumber\\
b_{1}=1,\,b_{n}=-1,\,b_{i}=0,\text{ for }i\neq1,n\nonumber\\
0\leq x_{ij}\leq1\nonumber
\end{gather}
\emph{Proof}: Given Lemmas 1 and 2, we know that the optimal solution
satisfies $\sum_{i}\sum_{j}r_{i}^{\ast}x_{ij}^{\ast}=\sum_{i}\sum_{j}%
e_{ij}x_{ij}^{\ast}-E_{1}$. Consequently, we can change the objective
(\ref{obj}) to the form below without affecting optimality:
\[
\min_{_{x_{ij,\text{ }}i,j\in\mathcal{N}}}\quad\sum_{i=1}^{n}\sum_{j=1}%
^{n}(\tau_{ij}+e_{ij}g)x_{ij}-E_{1}g
\]
Since $r_{i}$ no longer appears in the objective function and is only
contained in the energy dynamics (\ref{EiConv}), we can choose any $r_{i}$
satisfying the constraints (\ref{EiConv})-(\ref{Ei}) without affecting the
optimal objective function value. Therefore, $x_{ij}^{\ast}$ can be determined
by the following problem:
\begin{gather}
\min_{_{x_{ij,\text{ }}i,j\in\mathcal{N}}}\quad\sum_{i=1}^{n}\sum_{j=1}%
^{n}(\tau_{ij}+e_{ij}g)x_{ij}-E_{1}g\nonumber\\
s.t.\quad\sum_{j\in O(i)}x_{ij}-\sum_{j\in I(i)}x_{ji}=b_{i},\quad\text{for
each }i\in\mathcal{N}\nonumber\\
b_{1}=1,\,b_{n}=-1,\,b_{i}=0,\text{ for }i\neq1,n\nonumber\\
x_{ij}\in\{0,1\}\nonumber
\end{gather}
which is a typical shortest path problem formulation. Moreover, according to
the property of minimum cost flow problems [\cite{Hillier}], the above integer
programming problem is equivalent to the linear programming problem with the
integer restriction of $x_{ij}$ relaxed. Finally, since $E_{1}$ is given, the
problem reduces to (\ref{obj3}), which proves the theorem. $\blacksquare$

\subsection{Determination of optimal recharging amounts $r_{i}^{\ast}$}

Once we determine the optimal route, $P$, in (\ref{obj3}), it is relatively
easy to find a feasible solution for $r_{i}$, $i\in P$, to satisfy the
constraint (\ref{EiConv}), which is obviously non-unique in general. Then, we
can introduce a second objective into the problem, i.e., the minimization of
charging costs on the selected path, since charging prices normally vary over
stations. As before, we re-index nodes and define $P=\{1,...,n\}$. We denote
the charging price at node $i$ by $p_{i}$. Once an optimal route is
determined, we seek to control the energy recharging amounts $r_{i}$ to
minimize the total charging cost dependent on $p_{i}$, $i\in\mathcal{N}%
/\{n\}$. This can be formulated as a multistage optimal control problem:
\begin{gather}
\min_{r_{i},\text{ }i\in P}\quad\sum_{i\in P}p_{i}r_{i}\label{obj4}\\
s.t.\quad E_{i+1}=E_{i}+r_{i}-e_{i,i+1}\nonumber\\
0\leq E_{i}\leq B,\quad E_{1}\text{ given}\nonumber\\
r_{i}\geq0\text{ for all }i\in\mathcal{N}\nonumber
\end{gather}
This is a simple two-point boundary-value problem and can be easily solved by
discrete-time optimal control approaches [\cite{Ho}] or treating it as a
linear programming problem where $E_{i}$ and $r_{i}$ are both decision
variables. Due to space limitations, we omit numerical results providing
example solutions of the simple linear programming problem (\ref{obj3}) and
subsequent solutions of (\ref{obj4}).

Finally, we note that Theorem 1 holds under the assumption that charging nodes
are homogeneous in terms of charging speeds (i.e., the charging rate $1/g$ is
fixed). However, our analysis allows for inhomogeneous charging prices. The
case of node-dependent charging rates is the subject of ongoing work and can
be shown to still allow a decomposition of the MINLP, although we can no
longer generally obtain a LP.

\section{Multiple Vehicle Routing}

\label{sec3}

The results obtained for the single vehicle routing problem pave the way for
the investigation of multi-vehicle routing, where we seek to optimize a
system-wide objective by routing vehicles through the same network topology.
The main technical difficulty in this case is that we need to consider the
influence of traffic congestion on both traveling time and energy consumption.
A second difficulty is that of implementing an optimal routing policy. In the
case of a centrally controlled system consisting of mobile robots, sensors or
any type of autonomous vehicles this can be accomplished through appropriately
communicated commands. In the case of vehicles with individual drivers,
implementation requires signaling mechanisms and possibly incentive structures
to enforce desired routes assigned to vehicles, bringing up a number of
additional research issues. In the sequel, we limit ourselves to resolving the
first difficulty before addressing implementation challenges.

If we proceed as in the single vehicle case, i.e., determining a path
selection through $x_{ij}^{k}$, $i,j\in\mathcal{N}$, and recharging amounts
$r_{i}^{k}$, $i\in\mathcal{N}/\{n\}$ for all vehicles $k=1,\ldots,K$, for some
$K$, then the dimensionality of the solution space is prohibitive. Moreover,
the inclusion of traffic congestion effects introduces additional
nonlinearities in the dependence of the travel time $\tau_{ij}$ and energy
consumption $e_{ij}$ on the traffic flow through arc $(i,j)$, which now depend
on $x_{ij}^{1},\cdots,x_{ij}^{K}$. Instead, we will proceed by grouping
subsets of vehicles into $N$ \textquotedblleft subflows\textquotedblright%
\ where $N$ may be selected to render the problem manageable.

Let all vehicles enter the network at the origin node 1 and let $R$ denote the
rate of vehicles arriving at this node. Viewing vehicles as defining a
\emph{flow}, we divide them into $N$ \emph{subflows} (we will discuss the
effect of $N$ in Section 3.3), each of which may be selected so as to include
the same type of homogeneous vehicles (e.g., large vehicles vs smaller ones or
vehicles with the same initial energy). Thus, all vehicles in the same subflow
follow the same routing and recharging decisions so that we only consider
energy recharging at the subflow level rather than individual vehicles. Note
that asymptotically, as $N\rightarrow\infty$, we can recover routing at the
individual vehicle level.

Clearly, not all vehicles in our system are BPVs and are, therefore, not part
of our optimization process. These can be treated as uncontrollable
interfering traffic for our purposes and can be readily accommodated in our
analysis, as long as their flow rates are known. However, for simplicity, we
will assume here that every arriving vehicle is a BPV and joins a subflow.

Our objective is to determine optimal routes and energy recharging amounts for
each subflow of vehicles so as to minimize the total elapsed time of these
vehicle flows traveling from the origin to the destination. The decision
variables consist of $x_{ij}^{k}\in\{0,1\}$ for all arcs $(i,j)$ and subflows
$k=1,\ldots,N$, as well as charging amounts $r_{i}^{k}$ for all nodes
$i=1,\ldots,n-1$ and $k=1,\ldots,N$. Given traffic congestion effects, the
time and energy consumption on each arc depends on the values of $x_{ij}^{k}$
and the fraction of the total flow rate $R$ associated with each subflow $k$;
the simplest such flow allocation is one where each subflow is assigned $R/N$.
Let $\mathbf{x_{ij}}=(x_{ij}^{1},\cdots,x_{ij}^{N})^{T}$ and $\mathbf{r_{i}%
}=(r_{i}^{1},\cdots,r_{i}^{N})^{T}$. Then, we denote the traveling time and
corresponding energy consumption of the $k$th vehicle subflow on arc $(i,j)$
by $\tau_{ij}^{k}(\mathbf{x_{ij}})$ and $e_{ij}^{k}(\mathbf{x_{ij}})$
respectively. As already mentioned, $\tau_{ij}^{k}(\mathbf{x_{ij}})$ and
$e_{ij}^{k}(\mathbf{x_{ij}})$ can also incorporate the influence of
uncontrollable (non-BPV) vehicle flows, which can be treated as parameters in
these functions. Similar to the single vehicle case, we use $E_{i}^{k}$ to
represent the residual energy of subflow $k$ at node $i$, given by the
aggregated residual energy of all vehicles in the subflow. If the subflow does
not go through node $i$, then $E_{i}^{k}=0$. The problem formulation is as
follows:
\begin{gather}
\min_{\mathbf{x_{ij}},\mathbf{r_{i}},\,\,i,j\in\mathcal{N}}\quad\sum_{i=1}%
^{n}\sum_{j=1}^{n}\sum_{k=1}^{N}\left(  \tau_{ij}^{k}(\mathbf{x_{ij}}%
)+r_{i}^{k}gx_{ij}^{k}\right) \label{objM}\\
s.t.\text{ for each }k\in\{1,\dots,N\}:\nonumber\\
\quad\sum_{j\in O(i)}x_{ij}^{k}-\sum_{j\in I(i)}x_{ji}^{k}=b_{i}%
,\quad\text{for each }i\in\mathcal{N}\label{flowConvM}\\
b_{1}=1,\,b_{n}=-1,\,b_{i}=0,\text{ for }i\neq1,n\label{biM}\\
E_{j}^{k}=\sum_{i\in I(j)}(E_{i}^{k}+r_{i}^{k}-e_{ij}^{k}(\mathbf{x_{ij}%
}))x_{ij}^{k},\quad j=2,\dots,n\label{EiConvM}\\
E_{1}^{k}\text{ is given},\quad E_{i}^{k}\geq0,\quad\text{ for each }%
i\in\mathcal{N}\label{EiM}\\
x_{ij}^{k}\in\{0,1\},\quad r_{i}^{k}\geq0 \label{controlsM}%
\end{gather}
Obviously, this MINLP problem is difficult to solve. However, as in the
single-vehicle case, we are able to establish some properties that will allow
us to simplify it.

\subsection{Properties}

Even though the term $\tau_{ij}^{k}(\mathbf{x_{ij}})$ in the objective
function is no longer linear in general, for each subflow $k$ the constraints
(\ref{flowConvM})-(\ref{controlsM}) are still similar to the single-vehicle
case. Consequently, we can derive similar useful properties for this problem
in the form of the following two lemmas.

\textbf{Lemma 3: } For each subflow $k=1,\ldots,N$,
\begin{equation}
\sum_{i=1}^{n}\sum_{j=1}^{n}(r_{i}^{k}-e_{ij}^{k}(\mathbf{x_{ij}}))x_{ij}%
^{k}=E_{n}^{k}-E_{1}^{k} \label{Enk-E1k}%
\end{equation}

\textbf{Lemma 4:} If $\sum_{i=1}^{n}r_{i}^{k\ast}>0$ in the optimal routing
policy, then $E_{n}^{k\ast}=0$ for all $k=1,\ldots,N$.

The proofs of the above two lemmas are almost identical to those of Lemmas 1
and 2 respectively and are omitted. The only difference is that here the
analysis is focused on each vehicle \emph{subflow} instead of an
\emph{individual} vehicle. In view of Lemma 3, we can replace $\sum_{i=1}%
^{n}\sum_{j=1}^{n}r_{i}^{k}gx_{ij}^{k}$ in (\ref{objM}) by $(E_{n}^{k}%
-E_{1}^{k})g+\sum_{i=1}^{n}\sum_{j=1}^{n}e_{ij}^{k}(\mathbf{x_{ij}})gx_{ij}$
and eliminate, for all $k=1,\ldots,N$, the presence of $r_{i}^{k}$,
$i=1,\ldots,n-1$, from the objective function similar to the single-vehicle
case. Since $E_{1}^{k}$ is given, this leaves only the task of determining the
value of $E_{n}^{k}$. There are two possible cases: $(i)$ $\sum_{i}%
r_{i}^{k\ast}>0$, i.e., the $k$th vehicle subflow has to get recharged at
least once, and $(ii)$ $\sum_{i}r_{i}^{k\ast}=0$, i.e., $r_{i}^{k\ast}=0$ for
all $i$ and the $k$th vehicle subflow has adequate energy to reach the
destination without recharging.

Similar to the derivation of (\ref{leE1}), Case $(ii)$ results in a new
constraint $\sum_{i}\sum_{j}e_{ij}^{k}(\mathbf{x_{ij}})x_{ij}^{k}\leq
E_{1}^{k}$ for subflow $k$. However, since $e_{ij}^{k}(\mathbf{x_{ij}})$ now
depends on all $x_{ij}^{1},\ldots,x_{ij}^{N}$, the problem (\ref{objM}%
)-(\ref{controlsM}) with all $r_{i}^{k}=0$ is not as simple to solve as was
the case with (\ref{obj2})-(\ref{controls2}). Let us instead concentrate on
the more interesting Case $(i)$ for which Lemma 4 applies and we have
$E_{n}^{k\ast}=0$. Therefore, along with Lemma 3, we have for each
$k=1,\ldots,N$:
\[
\sum_{i=1}^{n}\sum_{j=1}^{n}r_{i}^{k}x_{ij}^{k}=\sum_{i=1}^{n}\sum_{j=1}%
^{n}e_{ij}^{k}(\mathbf{x_{ij}})x_{ij}^{k}-E_{1}^{k}%
\]
Then, proceeding as in Theorem 1, we can replace the original objective
function (\ref{objM}) and have the following new problem formulation to
determine $x_{ij}^{k\ast}$ for all $i,j\in\mathcal{N}$ and $k=1,\ldots,N$:
\begin{gather}
\min_{\mathbf{x_{ij}},\,\,i,j\in\mathcal{N}}\quad\sum_{i=1}^{n}\sum_{j=1}%
^{n}\sum_{k=1}^{N}\left(  \tau_{ij}^{k}(\mathbf{x_{ij}})+e_{ij}^{k}%
(\mathbf{x_{ij}})gx_{ij}^{k}\right) \label{objM2}\\
s.t.\text{ for each }k\in\{1,\dots,N\}:\nonumber\\
\quad\sum_{j\in O(i)}x_{ij}^{k}-\sum_{j\in I(i)}x_{ji}^{k}=b_{i}%
,\quad\text{for each }i\in\mathcal{N}\nonumber\\
b_{1}=1,\,b_{n}=-1,\,b_{i}=0,\text{ for }i\neq1,n\nonumber\\
x_{ij}^{k}\in\{0,1\}\nonumber
\end{gather}
Since the objective function is no longer necessarily linear in $x_{ij}^{k}$,
(\ref{objM2}) cannot be further simplified into an LP problem as in Theorem 1.
The computational effort required to Solve this problem heavily depends on the
dimensionality of the network and the number of subflows. Nonetheless, from
the transformed formulation above, we are still able to separate the
determination of routing variables $x_{ij}^{k}$ from recharging amounts
$r_{i}^{k}$. Similar to the single-vehicle case, once the routes are
determined, we can obtain any $r_{i}^{k}$ satisfying the energy constraints
(\ref{EiConvM})-(\ref{EiM}) such that $E_{n}^{k}=0$, thus preserving the
optimality of the objective value. To further determine $r_{i}^{k\ast}$, we
can introduce a second level optimization problem similar to the
single-vehicle case in (\ref{obj4}). Next, we will present an alternative
formulation for the original problem (\ref{objM})-(\ref{controlsM}) which
leads to a computationally simpler solution approach.

\subsection{Flow control formulation}

We begin by relaxing the binary variables in (\ref{controlsM}) by letting
$0\leq x_{ij}^{k}\leq1$. Thus, we switch our attention from determining a
single path for any subflow $k$ to several possible paths by treating
$x_{ij}^{k}$ as the normalized vehicle flow on arc $(i,j)$ for the $k$th
subflow. This is in line with many network routing algorithms in which
fractions $x_{ij}$ of entities are routed from a node $i$ to a neighboring
node $j$ using appropriate schemes ensuring that, in the long term, the
fraction of entities routed on $(i,j)$ is indeed $x_{ij}$. Following this
relaxation, the objective function in (\ref{objM}) is changed to:
\[
\min_{\mathbf{x_{ij}},\mathbf{r_{i}},\,\,i,j\in\mathcal{N}}\quad\sum_{i=1}%
^{n}\sum_{j=1}^{n}\sum_{k=1}^{N}\tau_{ij}^{k}(\mathbf{x_{ij}})+\sum_{i=1}%
^{n}\sum_{k=1}^{N}r_{i}^{k}g
\]
Moreover, the energy constraint (\ref{EiConvM}) needs to be adjusted
accordingly. Let $E_{ij}^{k}$ represent the fraction of residual energy of
subflow $k$ associated with the $x_{ij}^{k}$ portion of the vehicle flow
exiting node $i$. Therefore, the constraint (\ref{EiM}) becomes $E_{ij}%
^{k}\geq0$. We can now capture the relationship between the energy associated
with subflow $k$ and the vehicle flow as follows:
\begin{gather}
\left[  \sum_{h\in I(i)}(E_{hi}^{k}-e_{hi}^{k}(\mathbf{x_{ij}}))+r_{i}%
^{k}\right]  \cdot\frac{x_{ij}^{k}}{\sum_{h\in I(i)}x_{hi}^{k}}=E_{ij}%
^{k}\label{EiConvM2}\\
\frac{E_{ij}^{k}}{\sum_{j\in O(i)}E_{ij}^{k}}=\frac{x_{ij}^{k}}{\sum_{j\in
O(i)}x_{ij}^{k}} \label{same}%
\end{gather}
In (\ref{EiConvM2}), the energy values of different vehicle flows entering
node $i$ are aggregated and the energy corresponding to each portion exiting a
node, $E_{ij}^{k} $ , $j\in O(i)$, is proportional to the corresponding
fraction of vehicle flows, as expressed in (\ref{same}). Clearly, this
aggregation of energy leads to an approximation, since one specific vehicle
flow may need to be recharged in order to reach the next node in its path,
whereas another might have enough energy without being recharged. This
approximation foregoes controlling recharging amounts at the individual
vehicle level and leads to approximate solutions of the original problem
(\ref{objM})-(\ref{controlsM}). Several numerically based comparisons are
provided in the next section showing little or no loss of optimality relative
to the solution of (\ref{objM}).

Adopting this formulation with $x_{ij}^{k}\in\lbrack0,1]$ instead of
$x_{ij}^{k}\in\{0,1\}$, we obtain the following simpler nonlinear programming
problem (NLP):
\begin{gather}
\min_{\mathbf{x_{ij}},\mathbf{r_{i}},\,\,i,j\in\mathcal{N}}\quad\sum_{i=1}%
^{n}\sum_{j=1}^{n}\sum_{k=1}^{N}\tau_{ij}^{k}(\mathbf{x_{ij}})+\sum_{i=1}%
^{n}\sum_{k=1}^{N}r_{i}^{k}g\label{objM3}\\
s.t.\text{ for each }k\in\{1,\dots,N\}:\nonumber\\
\quad\sum_{j\in O(i)}x_{ij}^{k}-\sum_{j\in I(i)}x_{ji}^{k}=b_{i}%
,\quad\text{for each }i\in\mathcal{N}\label{flowConvM3}\\
b_{1}=1,\,b_{n}=-1,\,b_{i}=0,\text{ for }i\neq1,n\nonumber\\
\left[  \sum_{h\in I(i)}(E_{hi}^{k}-e_{hi}^{k}(\mathbf{x_{ij}}))+r_{i}%
^{k}\right]  \cdot\frac{x_{ij}^{k}}{\sum_{h\in I(i)}x_{hi}^{k}}=E_{ij}%
^{k}\label{EiConvM3}\\
\frac{E_{ij}^{k}}{\sum_{j\in O(i)}E_{ij}^{k}}=\frac{x_{ij}^{k}}{\sum_{j\in
O(i)}x_{ij}^{k}}\label{sameM3}\\
E_{ij}^{k}\geq0,\label{EiM3}\\
0\leq x_{ij}^{k}\leq1,\quad r_{i}^{k}\geq0 \label{controlsM3}%
\end{gather}
As in our previous analysis, we are able to eliminate $\mathbf{r_{i}}$ from
the objective function in (\ref{objM3}) as follows.

\textbf{Lemma 5: } For each subflow $k=1,\ldots,N$,
\begin{equation}
\sum_{i=1}^{n}r_{i}^{k}=\sum_{i=1}^{n}\sum_{j=1}^{n}e_{ij}^{k}(\mathbf{x_{ij}%
})+\sum_{i\in I(n)}E_{in}^{k}-\sum_{i\in O(1)}E_{1i}^{k}\nonumber
\end{equation}
\emph{Proof}: Summing (\ref{EiConvM3}) over all $i=1,\ldots,n$ gives%
\begin{align}
\sum_{i=1}^{n}r_{i}^{k}  &  =\sum_{i=1}^{n}\sum_{j=1}^{n}e_{ij}^{k}%
(\mathbf{x_{ij}})+\sum_{i=1}^{n}\sum_{j\in O(i)}E_{ij}^{k}\nonumber\\
&  \quad-\sum_{i=1}^{n}\sum_{h\in I(i)}E_{hi}^{k}\nonumber
\end{align}
and using (\ref{flowConvM3}),(\ref{sameM3}), we get
\begin{equation}
\sum_{i=1}^{n}r_{i}^{k}=\sum_{i=1}^{n}\sum_{j=1}^{n}e_{ij}^{k}(\mathbf{x_{ij}%
})+\sum_{i\in I(n)}E_{in}^{k}-\sum_{i\in O(1)}E_{1i}^{k}\nonumber
\end{equation}
which proves the lemma. $\blacksquare$

Similar to Lemma 3, we can easily see that if $\sum_{i}r_{i}^{k\ast}>0$ under
an optimal routing policy, then $\sum_{i\in I(n)}E_{in}^{k\ast}=0$. In
addition, $\sum_{i\in O(1)}E_{1i}^{k}=E_{1}^{k}$, which is given. We can now
transform the objective function (\ref{objM3}) into (\ref{objM4}) and
determine the optimal routes $x_{ij}^{k\ast}$ by solving the following NLP:
\begin{gather}
\min_{\substack{\mathbf{x_{ij}}\\i,j\in\mathcal{N}}}\quad\sum_{k=1}^{N}\left(
\sum_{i=1}^{n}\sum_{j=1}^{n}\left[  \tau_{ij}^{k}(\mathbf{x_{ij}})+e_{ij}%
^{k}(\mathbf{x_{ij}})g\right]  -E_{1}^{k}\right) \label{objM4}\\
s.t.\text{ for each }k\in\{1,\dots,N\}:\nonumber\\
\quad\sum_{j\in O(i)}x_{ij}^{k}-\sum_{j\in I(i)}x_{ji}^{k}=b_{i}%
,\quad\text{for each }i\in\mathcal{N}\nonumber\\
b_{1}=1,\,b_{n}=-1,\,b_{i}=0,\text{ for }i\neq1,n\nonumber\\
0\leq x_{ij}^{k}\leq1\nonumber
\end{gather}
The values of $r_{i}^{k}$, $i=1,\ldots,n$, $k=1,\ldots,N$, can be determined
so as to satisfy the energy constraints (\ref{EiConvM3})-(\ref{EiM3}), and
they are obviously not unique. We may then proceed with a second-level
optimization problem to determine optimal values similar to Section 2.2.

\subsection{Numerical Examples}

\begin{figure}[b]
\begin{center}
\includegraphics[scale=0.6]{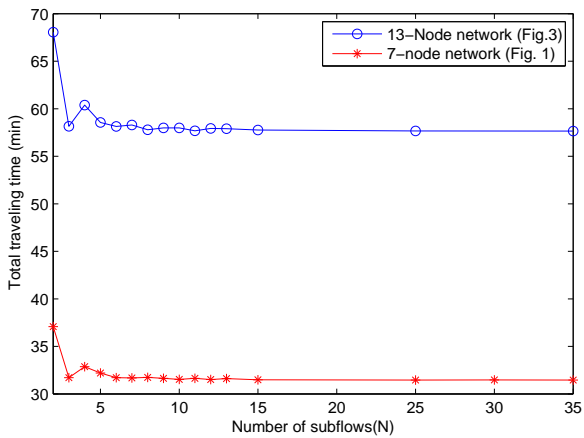}
\end{center}
\caption{Performance as a function of $N$ (No. of subflows)}%
\label{NvsObj}%
\end{figure}

We consider a specific example which includes traffic congestion and energy
consumption functions. The relationship between the speed and density of a
vehicle flow is typically estimated as follows (see \cite{Ioannou96}):
\begin{equation}
v(k(t))=v_{f}\bigg(1-\left(  \frac{k(t)}{k_{jam}}\right)  ^{p}\bigg)^{q}
\label{VvsD}%
\end{equation}
where $v_{f}$ is the reference speed on the road without traffic, $k(t)$
represents the density of vehicles on the road at time $t$ and $k_{jam}$ the
saturated density for a traffic jam. The parameters $p$ and $q$ are
empirically identified for actual traffic flows. In our multi-vehicle routing
problem, we are interested in the relationship between the density of the
vehicle flow and traveling time on an arc $(i,j)$, i.e., $\tau_{ij}%
^{k}(\mathbf{x_{ij}})$. Given a network topology (i.e., a road map), the
distances $d_{ij}$ between nodes are known. Moreover, we do not include
uncontrollable vehicle flows in our example for simplicity. In our approach,
we need to identify $N$ subflows and we do so by evenly dividing the entire
vehicle inflow into $N$ subflows, each of which has $R/N$ vehicles per unit
time. Thus, $k_{jam}$ in this case can be set as $N$, implying that we do not
want all vehicles to go through the same path, hence the the arc $(i,j)$
density is $\sum_{k}x_{ij}^{k}$. Therefore, the time subflow $k$ spends on arc
$(i,j)$ becomes
\[
\tau_{ij}^{k}(\mathbf{x_{ij}})=\frac{d_{ij}\cdot x_{ij}^{k}\cdot\frac{R}{N}%
}{v_{f}(1-(\frac{\sum_{k}x_{ij}^{k}}{N})^{p})^{q}}%
\]
As for $e_{ij}^{k}(\mathbf{x_{ij}})$, we assume the energy consumption rates
of subflows on arc $(i,j)$ are all identical, proportional to the distance
between nodes $i$ and $j$, giving
\[
e_{ij}^{k}(\mathbf{x_{ij}})=e\cdot d_{ij}\cdot\frac{R}{N}%
\]
Therefore, we aim to solve the multi-vehicle routing problem using
(\ref{objM2}) which in this case becomes:
\begin{gather}
\min_{\substack{x_{ij}^{k}\\i,j\in\mathcal{N}}}\quad\sum_{i=1}^{n}\sum
_{j=1}^{n}\sum_{k=1}^{N}\left(  \frac{d_{ij}x_{ij}^{k}\frac{R}{N}}%
{v_{f}(1-(\frac{\sum_{k}x_{ij}^{k}}{N})^{p})^{q}}+egd_{ij}\frac{R}{N}%
x_{ij}^{k}\right) \label{objSp}\\
s.t.\text{ for each }k\in\{1,\dots,N\}:\nonumber\\
\quad\sum_{j\in O(i)}x_{ij}^{k}-\sum_{j\in I(i)}x_{ji}^{k}=b_{i}%
,\quad\text{for each }i\in\mathcal{N}\nonumber\\
b_{1}=1,\,b_{n}=-1,\,b_{i}=0,\text{ for }i\neq1,n\nonumber\\
x_{ij}^{k}\in\{0,1\}\nonumber
\end{gather}
For simplicity, we let $v_{f}=1$ mile/min, $R=1$ vehicle/min, $p=2,\,\,q=2$
and $e\cdot g=1$. The network topology used is that of Fig.\ref{SampleNet},
where the distance of each arc is shown in Tab. \ref{table1}.
\begin{table}[pt]
\caption{$d_{ij}$ values for network of Fig. \ref{SampleNet} ($miles$)}%
\label{table1}%
\centering
\begin{tabular}
[c]{|c|c|c|c|c|c|c|c|c|c|}\hline\hline
$d_{12}$ & $d_{14}$ & $d_{15}$ & $d_{23}$ & $d_{24}$ & $d_{46}$ & $d_{56}$ &
$d_{37}$ & $d_{47}$ & $d_{67}$\\\hline
5 & 6.2 & 7 & 3.5 & 5 & 3.6 & 4.3 & 6 & 6 & 4\\\hline
\end{tabular}
\end{table}To solve the nonlinear binary programming problem (\ref{objSp}), we
use the optimization solver \emph{Opti} (MATLAB toolbox for optimization). The
results are shown in Tab. \ref{table2} for different values of $N=1,\ldots
,30$. As shown in Tab. \ref{table2}, vehicles are mainly distributed through
three routes and the traffic congestion effect makes the flow distribution
differ from following the shortest path. The number of decision variables
(hence, the solution search space) rapidly increases with the number of
subflows. However, looking at Fig. \ref{NvsObj} which gives the performance in
terms of our objective function in (\ref{objSp}) as a function of the number
of subflows, observe that the optimal objective value quickly converges around
$N=10$. Thus, even though the best solution is found when $N=25$, a
near-optimal solution can be determined under a small number of subflows. This
suggests that one can rapidly approximate the asymptotic solution of the
multi-vehicle problem (dealing with individual vehicles routed so as to
optimize a systemwide objective) based on a relatively small value of $N$.
\begin{table}[h]
\caption{Numerical results for sample problem}%
\label{table2}
\begin{center}%
\begin{tabular}
[c]{c|c|c}\hline\hline
N & 1 & 2\\\hline
obj & 1.22e9 & 37.077\\\hline
routes & $1\rightarrow4\rightarrow7$ & $%
\begin{array}
[c]{c}%
1\rightarrow4\rightarrow7\\
1\rightarrow2\rightarrow3\rightarrow7
\end{array}
$\\\hline\hline
N & 3 & 4\\\hline
obj & 31.7148 & 32.8662\\\hline
routes & $%
\begin{array}
[c]{c}%
(1\rightarrow4\rightarrow7)\\
1\rightarrow2\rightarrow3\rightarrow7\\
1\rightarrow5\rightarrow6\rightarrow7
\end{array}
$ & $%
\begin{array}
[c]{c}%
(1\rightarrow4\rightarrow7)\times2\\
1\rightarrow2\rightarrow3\rightarrow7\\
1\rightarrow5\rightarrow6\rightarrow7
\end{array}
$\\\hline\hline
N & 5 & 6\\\hline
obj & 32.1921 & 31.7148\\\hline
routes & $%
\begin{array}
[c]{c}%
(1\rightarrow4\rightarrow7)\times2\\
(1\rightarrow2\rightarrow3\rightarrow7)\times2\\
1\rightarrow5\rightarrow6\rightarrow7
\end{array}
$ & $%
\begin{array}
[c]{c}%
(1\rightarrow4\rightarrow7)\times2\\
(1\rightarrow2\rightarrow3\rightarrow7)\times2\\
(1\rightarrow5\rightarrow6\rightarrow7)\times2
\end{array}
$\\\hline\hline
N & 10 & 15\\\hline
obj & 31.5279 & 31.4851\\\hline
routes & $%
\begin{array}
[c]{c}%
(1\rightarrow4\rightarrow7)\times4\\
(1\rightarrow2\rightarrow3\rightarrow7)\times3\\
(1\rightarrow5\rightarrow6\rightarrow7)\times3
\end{array}
$ & $%
\begin{array}
[c]{c}%
(1\rightarrow4\rightarrow7)\times5\\
(1\rightarrow2\rightarrow3\rightarrow7)\times5\\
(1\rightarrow5\rightarrow6\rightarrow7)\times4\\
(1\rightarrow4\rightarrow6\rightarrow7)\times1
\end{array}
$\\\hline\hline
N & 25 & 30\\\hline
obj & 31.4513 & 31.4768\\\hline
routes & $%
\begin{array}
[c]{c}%
(1\rightarrow4\rightarrow7)\times9\\
(1\rightarrow2\rightarrow3\rightarrow7)\times8\\
(1\rightarrow5\rightarrow6\rightarrow7)\times7\\
(1\rightarrow4\rightarrow6\rightarrow7)\times1
\end{array}
$ & $%
\begin{array}
[c]{c}%
(1\rightarrow4\rightarrow7)\times11\\
(1\rightarrow2\rightarrow3\rightarrow7)\times10\\
(1\rightarrow5\rightarrow6\rightarrow7)\times8\\
(1\rightarrow4\rightarrow6\rightarrow7)\times1
\end{array}
$\\\hline\hline
\end{tabular}
\end{center}
\end{table}

Next, we obtain a solution to the same problem (\ref{objSp}) using the
alternative NLP formulation (\ref{objM4}) where $0\leq x_{ij}^{k}\leq1$. Since
in this example all subflows are identical, we can further combine all
$x_{ij}^{k}$ over each arc $(i,j)$, which leads to the following $N$-subflow
relaxed problem:
\begin{gather}
\min_{x_{ij},\,\,i,j\in\mathcal{N}}\quad\sum_{i=1}^{n}\sum_{j=1}^{n}\left(
\frac{d_{ij}x_{ij}R}{v_{f}(1-(x_{ij})^{p})^{q}}+egd_{ij}Rx_{ij}\right)
\label{objSp2}\\
s.t.\quad\sum_{j\in O(i)}x_{ij}-\sum_{j\in I(i)}x_{ji}=b_{i},\quad\text{for
each }i\in\mathcal{N}\nonumber\\
b_{1}=1,\,b_{n}=-1,\,b_{i}=0,\text{ for }i\neq1,n\nonumber\\
0\leq x_{ij}\leq1\nonumber
\end{gather}
This is a relatively easy to solve NLP problem. Using the same parameter
settings as before, we obtain the objective value of $31.4465$ mins and the
optimal routes are:
\begin{align}
35.88\%\text{ of vehicle flow: }  &  (1\rightarrow4\rightarrow7)\nonumber\\
31.74\%\text{ of vehicle flow: }  &  (1\rightarrow2\rightarrow3\rightarrow
7)\nonumber\\
27.98\%\text{ of vehicle flow: }  &  (1\rightarrow5\rightarrow6\rightarrow
7)\nonumber\\
4.44\%\text{ of vehicle flow: }  &  (1\rightarrow4\rightarrow6\rightarrow
7)\nonumber
\end{align}
Compared to the best solution ($N=25$) in Tab. \ref{table2} and Fig.
\ref{NvsObj}, the difference in objective values between the integer and
flow-based solutions is less than $0.1\%$. This supports the effectiveness of
a solution based on a limited number of subflows in the MINLP problem.

\textbf{Performance improvement over uncontrolled traffic systems}. Next, we
address the extent to which this optimization approach offers improvements
over an uncontrolled traffic network. We simulate the vehicle routing problem
on the discrete event simulator, MATLAB/SimEvents, where the vehicle arrivals
to the source are randomly generated with a random initial energy. As a simple
example, we model the routing for each vehicle at each node to be round-robin,
while the recharging amount of the vehicle is just adequate to reach the next
node. The objective value of such an uncontrolled routing policy for network shown in Fig.\ref{SampleNet} is $38.524$
mins, compared to our optimal policy which gave $31.451$ mins, an improvement
of $18.36\%$.

\begin{figure}[ptbh]
\begin{center}
\includegraphics[scale=0.5]{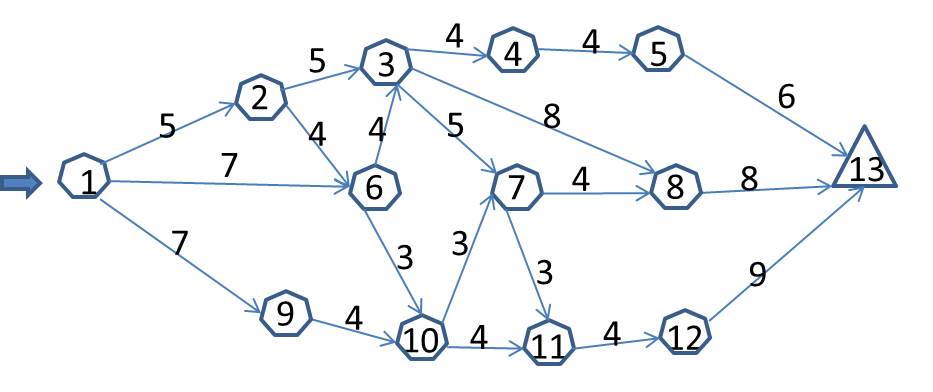}
\end{center}
\caption{A 13-node network example for routing with recharging nodes.}%
\label{bigGraph}%
\end{figure}

\textbf{Larger networks}. We have also considered a more topologically
complicated network with 13 nodes and 20 arcs as shown in Fig. \ref{bigGraph}.
The number on each arc indicates the distance between adjacent nodes. We
assume all other numerical values to be similar to the previous example. Fig.
\ref{NvsObj} shows the performance in terms of the objective function in
(\ref{objSp}) vs the number of subflows for this network. We can see that the
optimal objective value converges around $N=10$.

Now, let us solve the $N$-subflow relaxed problem (\ref{objSp2}) for this
network with the same parameter settings as before to check for its accuracy.
We obtain the optimal objective function value as 57.6326 which is almost
equal to the optimal traveling time of 57.6489
 obtained for $N=35$ in the
MINLP formulation. The optimal routing probabilities are as follows:\newline%

\begin{align*}
34.77\%\text{ of vehicle flow: }  &  (1\rightarrow2\rightarrow3\rightarrow4
\rightarrow5 \rightarrow13)\\
27.52\%\text{ of vehicle flow: }  &  (1\rightarrow9\rightarrow10\rightarrow11
\rightarrow12 \rightarrow13)\\
24.89\%\text{ of vehicle flow: }  &  (1\rightarrow6\rightarrow10\rightarrow7
\rightarrow8 \rightarrow13)\\
10.807\%\text{ of vehicle flow: }  &  (1\rightarrow6\rightarrow3\rightarrow8
\rightarrow13)\\
1.7\%\text{ of vehicle flow: } &  (1\rightarrow9\rightarrow10\rightarrow7
\rightarrow8 \rightarrow13)\\
0.313\%\text{ of vehicle flow: }  &  (1\rightarrow6\rightarrow3\rightarrow4
\rightarrow5 \rightarrow13)\\
\end{align*}

\textbf{CPU time Comparison}. Based on our simulation results we conclude that
the flow control formulation is a good approximation of the original MINLP
problem. Tab. \ref{table3} compares the computational effort in terms of CPU
time for both formulations to find optimal routes for the two sample networks
we have considered. Our results show that the flow control formulation results
in a reduction of about 5 orders of magnitude in CPU time with virtually identical objective function values.

\begin{table}[h]
\caption{CPU time for sample problems}%
\label{table3}
\begin{center}%
\begin{tabular}
[c]{c|c|c|c}\hline\hline
\textbf{Fig.\ref{SampleNet} Net.} & MINLP & MINLP & NLP approx.\\\hline
N & 2 & 10(near opt) & -\\\hline
obj & 37.083 & 31.5319 & 31.4504\\\hline
CPU time(sec) & 312 & 9705 & 0.07\\\hline\hline
\textbf{Fig.\ref{bigGraph} Net.} & MINLP & MINLP & NLP approx.\\\hline
N & 2 & 15(near opt) & -\\\hline
obj & 68.055 & 57.764 & 57.6326\\\hline
CPU time(sec) & 820 & 10037 & 0.2\\\hline\hline
\end{tabular}
\end{center}
\end{table}\textbf{Effect of recharging speed on optimal routes}. Once we
determine the optimal routes, we can also ascertain the total time spent
traveling and recharging respectively, i.e., the first and second terms in
(\ref{objSp2}). Obviously the value of $e\cdot g$, which captures the
recharging speed, determines the proportion of traveling and recharging amount
as well as the route selection. As shown in Tab. \ref{table4}, the larger the
product $e\cdot g$ is, the slower the recharging speed, therefore the more
weighted the recharging time in the objective function becomes. In this case,
flows tend to select shortest paths in terms of energy consumption.
Conversely, if the recharging speed is fast, the routes are selected to
prioritize the traveling time on paths. \begin{table*}[ht]
\caption{Numerical results for different values of $e\cdot g$ for network of
Fig. \ref{SampleNet}}%
\label{table4}%
\centering
\begin{tabular}
[c]{c|c|c|c}\hline\hline
$e\cdot g$ & 0.1 & 1 & 10\\\hline
total time & 18.9417 & 31.4465 & 154.4777\\\hline
time on paths & 17.5471 & 17.5791 & 19.4510\\\hline
time at stations & 1.3946 & 13.8674 & 135.0267\\\hline
optimal routes & $%
\begin{array}
[c]{cl}%
31.53\%: & (1\rightarrow2\rightarrow3\rightarrow7)\\
32.97\%: & (1\rightarrow4\rightarrow7)\\
28.58\%: & (1\rightarrow5\rightarrow6\rightarrow7)\\
5.78\%: & (1\rightarrow4\rightarrow6\rightarrow7)\\
1.14\%: & (1\rightarrow2\rightarrow4\rightarrow7)
\end{array}
$ & $%
\begin{array}
[c]{cl}%
31.74\%: & (1\rightarrow2\rightarrow3\rightarrow7)\\
35.88\%: & (1\rightarrow4\rightarrow7)\\
27.98\%: & (1\rightarrow5\rightarrow6\rightarrow7)\\
4.4\%: & (1\rightarrow4\rightarrow6\rightarrow7)
\end{array}
$ & $%
\begin{array}
[c]{cl}%
32.35\%: & (1\rightarrow2\rightarrow3\rightarrow7)\\
49.63\%: & (1\rightarrow4\rightarrow7)\\
18.02\%: & (1\rightarrow5\rightarrow6\rightarrow7)
\end{array}
$\\\hline\hline
\end{tabular}
\end{table*}

\section{Conclusions and future work}

\label{sec4}

We have introduced energy constraints into the vehicle routing problem, and
studied the problem of minimizing the total elapsed time for vehicles to reach
their destinations by determining routes as well as recharging amounts when
there is no adequate energy for the entire journey. For a single vehicle, we
have shown how to decompose this problem into two simpler problems. For a
multi-vehicle problem, where traffic congestion effects are considered, we
used a similar approach by aggregating vehicles into subflows and seeking
optimal routing decisions for each such subflow. We also developed an
alternative flow-based formulation which yields approximate solutions with a
computational cost reduction of several orders of magnitude, so they can be
used in problems of large dimensionality. Numerical examples show these
solutions to be near-optimal. We have also found that a low number of subflows
is adequate to obtain convergence to near-optimal solutions, making the
multi-subflow strategy particularly promising.

Our ongoing work introduces different characteristics into the charging
stations, such as recharging speeds and queueing capacities. In this case, we
can show that a similar decomposition still holds, although we can no longer
obtain an LP problem. We also believe that extensions to multiple vehicle
origins and destinations are straight-forward, as is the case where only a
subset of nodes has recharging resources or not all vehicles in the network are
BPVs. Finally, we are exploring extensions into stochastic vehicle flows which
can incorporate various random effects.

\bibliographystyle{IEEEtran}
\bibliography{IEEEabrv,renowang}

\end{document}